\documentclass{patmorin}
\usepackage{pat}
\usepackage{paralist}
\usepackage{dsfont}  % for \mathds{A}
\usepackage[utf8x]{inputenc}
\usepackage{skull}
\usepackage{paralist}
\usepackage{graphicx}
\usepackage[noend]{algorithmic}

\usepackage[normalem]{ulem}
\usepackage{cancel}
\usepackage{enumitem}

\usepackage{pifont}
\usepackage{todonotes}

\usepackage{relsize}

\usepackage[longnamesfirst,numbers,sort&compress]{natbib}

\usepackage[mathlines]{lineno}
\newcommand*\patchAmsMathEnvironmentForLineno[1]{%
 \expandafter\let\csname old#1\expandafter\endcsname\csname #1\endcsname
 \expandafter\let\csname oldend#1\expandafter\endcsname\csname end#1\endcsname
 \renewenvironment{#1}%
    {\linenomath\csname old#1\endcsname}%
    {\csname oldend#1\endcsname\endlinenomath}}%
\newcommand*\patchBothAmsMathEnvironmentsForLineno[1]{%
 \patchAmsMathEnvironmentForLineno{#1}%
 \patchAmsMathEnvironmentForLineno{#1*}}%
\AtBeginDocument{%
\patchBothAmsMathEnvironmentsForLineno{equation}%
\patchBothAmsMathEnvironmentsForLineno{align}%
\patchBothAmsMathEnvironmentsForLineno{flalign}%
\patchBothAmsMathEnvironmentsForLineno{alignat}%
\patchBothAmsMathEnvironmentsForLineno{gather}%
\patchBothAmsMathEnvironmentsForLineno{multline}%
}

\title{\MakeUppercase{Odd Colourings of Graph Products}\thanks{This research was partly funded by NSERC.}}
\author{%
  Vida Dujmović\thanks{Department of Computer Science and Electrical Engineering, University of Ottawa}\qquad
  Pat Morin\thanks{School of Computer Science, Carleton University}\qquad
  Saeed Odak\footnotemark[2]
}

\date{}

\usepackage{tabularx}

\begin{document}

% \begin{titlepage}
\maketitle

\begin{abstract}
  The odd colouring number is a new graph parameter introduced by \citet{petrusevski.skrekovski:colorings}.  In this note, we show that graphs with so called product structure have bounded odd-colouring number. By known results on the product structure of $k$-planar graphs, this implies that $k$-planar graphs have bounded odd-colouring number, which answers a question of \citet{cranston.lafferty.ea:note}.
\end{abstract}
% \end{titlepage}

% \pagenumbering{roman}
% \tableofcontents
%
% \newpage
% \pagenumbering{arabic}

\section{Introduction}

Let $G$ be a graph.  A (not necessarily proper)\footnote{$\varphi$ is a proper colouring of $G$ if $vw\in E(G)$ implies that $\varphi(v)\neq\varphi(w)$.} vertex colouring $\varphi:V(G)\to\N$ is \emph{odd} if the neighbourhood of each non-isolated vertex of $G$ contains a colour that occurs an odd number of times.  More precisely, if $N_G(v):=\{w\in V(G):vw\in E(G)\}$ denotes the neighbourhood of a vertex $v$ in $G$, then $\varphi$ is an odd colouring of $G$ if and only if, for each $v\in V(G)$ with $|N_G(v)|>0$, there exists a colour $\alpha$ such that $|\{w\in N_G(v): \varphi(w)=\alpha\}|$ is odd.

Odd colourings were recently introduced by \citet{petrusevski.skrekovski:colorings}, who showed that every planar graph $G$ has a proper odd colouring using at most $9$ colours\footnote{We say that a colouring $\varphi$ of $G$ \emph{uses $c$ colours} if $c=|\{\varphi(v):v\in V(G)\}$.}, and conjectured that $5$ colours always suffice.  \citet{caro.petrusevski.ea:remarks} showed that $5$ colours always suffice for outerplanar graphs and showed that $8$ colours always suffice for some special cases of planar graphs.  Building on the work of \citet{caro.petrusevski.ea:remarks}, \citet{petr.portier:odd} showed that every planar graph has an odd colouring using at most $8$ colours.

A minor-closed family of graphs\footnote{A graph $M$ is a \textit{minor} of a graph $G$ if a graph isomorphic to $M$ can be obtained from a subgraph of $G$ by contracting edges. A class $\mathcal{G}$ of graphs is \emph{minor-closed} if for every graph $G\in\mathcal{G}$, every minor of $G$ is in $\mathcal{G}$.} $\mathcal{G}$ is \emph{$d$-degenerate} if every graph in $\mathcal{G}$ contains a vertex of degree at most $d$.
% (Note that any proper minor-closed family of graphs is $d$-degenerate for some fixed $d$.) A vertex of degree zero is called an \emph{isolated} vertex.
\citet{cranston.lafferty.ea:note} proved that any graph from a $d$-degenerate minor-closed family of graphs has a proper odd colouring using at most $2d+1$ colours.  This result, which has a short and elegant proof, includes outerplanar graphs and, more generally, partial $2$-trees (with $d=2$); planar graphs (with $d=5$); graphs embeddable on surfaces of Euler genus $g$; and graphs of treewidth at most $t$ (with $d=t$).

\citet{cranston.lafferty.ea:note} also consider $1$-planar graphs, which do not form a minor-closed family, and show that any $1$-planar graph\footnote{A graph is $k$-planar if it has an embedding in the plane in which no edge contains a vertex other than its endpoints and each edge is involved in at most $k$ crossings with other edges.} has a proper odd colouring using at most $31$ colours.  They ask if this can be extended to $k$-planar graphs for general $k>1$.  Our main result, \cref{new_result} below, is a more general results for graphs having \emph{product structure} and implies that any $k$-planar graph has a proper odd colouring using $O(k^5)$ colours.
%
% \begin{thm}\label{old_result}
%   Let $\mathcal{H}$ be a $d$-degenerate minor closed family of graphs and let $H$ be a member of $\mathcal{H}$ with no isolated vertices.  Then $H$ has a proper odd colouring $\varphi:V(H)\to\{1,\ldots,2d+1\}$.
% \end{thm}

For two graphs $A$ and $B$, the \emph{strong product} $A\boxtimes B$ of $A$ and $B$ is the graph with vertex set $V(A\boxtimes B):=V(A)\times V(B)$ and that contains an edge with endpoints $(x_1,y_1)$ and $(x_2,y_2)$ if and only if
\begin{inparaenum}[(i)]
  \item $x_1x_2\in E(A)$ and $y_1=y_2$;
  \item $x_1=x_2$ and $y_1y_2\in E(B)$; or
  \item $x_1x_2\in E(A)$ and $y_1y_2\in E(B)$.
\end{inparaenum}
A \emph{$t$-tree} $H$ is a graph that is either a clique on $t+1$ vertices or a graph that contains a vertex $v$ of degree $t$ whose neighbours form a clique and such that $H-\{v\}$ is a $t$-tree.  The following is the main result in this paper:

\begin{thm}\label{new_result}
  Let $H$ be a $t$-tree, let $P$ be a path, and let $G$ be a subgraph of $H\boxtimes P$. Then $G$ has a proper odd colouring of $G$ that uses at most $8t+4$ colours.
\end{thm}

A number of graph families are known to exhibit so called \emph{product structure} like that required of the graph $G$ in \cref{new_result}.  For example, for every planar graph $G$ there exists a $6$-tree $H$ and a path $P$ such that $G$ is isomorphic to a subgraph of $H\boxtimes P$ \cite{ueckerdt.wood.ea:improved}.  Similar results hold (with constants other than $6$) for graphs of bounded Euler genus, apex-minor free graphs, and bounded-degree graphs from proper-minor-closed families \cite{dujmovic.joret.ea:planar,dujmovic.esperet.ea:clustered}.  Most relevant for the current discussion is the following theorem of \citet{dujmovic.morin.ea:structure}:

\begin{thm}[\citet{dujmovic.morin.ea:structure}]\label{k_planar}
  For every $k$-planar graph $G$ there exists an $O(k^5)$-tree $H$ and a path $P$ such that $G$ is isomorphic to a subgraph of $H\boxtimes P$.
\end{thm}

Combining this with \cref{new_result}, we immediate obtain the following corollary, answering the question posed by \citet{cranston.lafferty.ea:note}.

\begin{cor}
  Every $k$-planar graph $G$ has a proper odd colouring that uses $O(k^5)$ colours.
\end{cor}

\section{Proof of \cref{new_result}}

\begin{proof}[Proof of \cref{new_result}]
  Let $y_1,\ldots,y_h$ be the vertices of $P$, in order.  To avoid a boring edge case, we extend $P$ by one vertex in each direction, so that the vertices $y_0$ and $y_{h+1}$ are defined.

  Let $x_1,\ldots,x_r$ be the vertices of $H$ ordered so that $x_1,\ldots,x_{t}$ is a clique and, for each $i\in\{t+1,\ldots,r\}$, the vertices in $C_{x_i}:=N_G(x_i)\cap\{x_1,\ldots,x_{i-1}\}$ form a clique of order $t$.  We make crucial use of the following well-known property of $t$-trees:\footnote{(\ding{74}) follows from the fact that $S:=C_{x_i}\cup\{x_i\}$ separates $\{x_1,\ldots,x_i\}\setminus S$ from $x_{i'}$.  In the language of rooted tree-decompositions, the node whose bag contains $S$ is an ancestor of all nodes whose bags contains $x_{i'}$.}
  \begin{compactitem}[(\ding{74})]
    \item If $1\le i \le j\le r$ and $x_i \in C_{x_{j}}$ then $C_{x_i}\cup\{x_i\}\supseteq C_{x_{j}}\cap\{x_1,\ldots,x_{i}\}$.
  \end{compactitem}

  With this setup out of the way, we can proceed with the proof, which is by induction on the number of vertices of $G$.  We will prove the following stronger statement:  There exists a proper odd colouring $\varphi:V(G)\to\{1,\ldots,8t+4\}$ of $G$ that satisifies the following additional condition:
  \begin{compactitem}[(\ding{96})]
    \item For each $(x_i,y_j)\in V(H\boxtimes P)$, define
    \[
      C_{(x_i,y_j)}:=V(G)\cap \left(\left(C_{x_i}\times\{y_{j-1}, y_{j},y_{j+1}\}\right)\cup\{(x_i,y_i),(x_i,y_{j-1})\}\right)
      \enspace .
    \]
    Then, for each $v\in V(H\boxtimes P)$ the vertices in $C_v$ receive distinct colours, i.e., $\varphi(u)\neq \varphi(w)$ for each distinct $u,w\in C_v$.
  \end{compactitem}
  Note that, for any edge $vw$ of $G$, $v,w \in C_w$ or $v,w\in C_v$.  Therefore, (\ding{96}) implies that the colouring $\varphi$ is a proper colouring of $G$.

  The base case, in which $G$ has no vertices, is trivial. Therefore, we may assume that $G$ has at least one vertex.  Let $(i,j)\in\{1,\ldots,r\}\times\{1,\ldots,h\}$ be the lexicographically largest pair such that $v:=(x_i,y_j)\in V(G)$.  Observe that $|N_G(v)|\le 3t+1$ since, by the maximality of $(i,j)$,  $N_G(v)\subseteq C_v\setminus\{v\}$.

  Consider the graph $G':= G-\{v\}$. Since $G'$ is a subgraph of $H\boxtimes P$ with fewer vertices than $G$, the inductive hypothesis implies that there exists an odd colouring $\varphi$ of $G'$ that satisfies (\ding{96}).  We will now extend $\varphi$ to a colouring of $G$ by listing colours that we may not choose for $\varphi(v)$:
  \begin{itemize}
    \item To guarantee that $\varphi$ satisfies (\ding{96}), observe that assigning a colour to $v$ can only violate (\ding{96}) if it does so for some vertex $(x_{i'},y_{j'})\in V(H\boxtimes P)$ with $v\in C_{(x_{i'},y_{j'})}$.  By (\ding{74}) and the maximality of $(i,j)$, if $v\in C_{(x_{i'},y_{j'})}$ then $|j-j'|\le 1$ and $C_{x_{i}}\cup\{x_i\}\supseteq C_{x_{i'}}$.  Therefore,
    \[
        C_{(x_{i'},y_{j'})}\setminus\{v\}\subseteq
          \left(C_{x_i} \times \{y_{j-2},y_{j-1},y_j,y_{j+1},y_{j+2}\}\right)
          \cup \{(x_i,y_{j-2}),(x_i,y_{j-1})\} =: R \enspace .
    \]
    Therefore, in order to satisfy (\ding{96}) it is sufficient to choose $\varphi(v)$ so that $\varphi(v)\neq\varphi(w)$ for each $w\in R$.  Let $X:=\{\varphi(w): w\in R\}$ and observe that $|X|\le |R| \le 5t+2$.

    Furthermore, if $v$ is not an isolated vertex of $G$ then (\ding{96}) ensures that some colour occurs exactly once in $N_G(v)$. Therefore, to ensure that $\varphi$ is an odd colouring of $G$, we need only choose some $\varphi(v)\not\in X$ in such a way that some colour appears an odd number of times in $N_G(w)$ for each $w\in N_G(v)$, which is what we do next.

    \item To guarantee that $\varphi$ is an odd colouring of $G$, consider each vertex $w\in N_{G}(v)$.  If there is exactly one colour $\alpha\in\{1,\ldots,8t+4\}$ that occurs an odd number of times in $N_{G'}(w)$, then define $Y_{w} := \{\alpha\}$; otherwise define $Y_{w}:=\emptyset$. Now let $Y:=\bigcup_{w\in N_{G}(v)} Y_{w}$ and observe that $|Y|\le |N_G(v)|\le 3t+1$.
    If we choose $\varphi(v)\not\in Y$ then , for each $w\in N_{G}(v)$ the following holds:
    \begin{itemize}
      \item If $Y_{w}=\{\alpha\}$ then $\varphi(v)\neq\alpha$. Therefore, the colour $\alpha$ appears an odd number of times in $N_{G}(w)$ since it appears an odd number of times in $N_{G'}(w)=N_G(w)\setminus\{v\}$.
      \item If $Y_{w}=\emptyset$ then either:
      \begin{itemize}
        \item No colour appears an odd number of times in $N_{G'}(w)$.  Therefore the colour $\varphi(v)$ appears an even number of times in $N_{G'}(w)$, so $\varphi(v)$ appears an odd number of times in $N_G(w)=N_{G'}(w)\cup\{v\}$.
        \item At least two colours $\alpha$ and $\beta$ each appear an odd number of times in $N_{G'}(w)$.  Therefore each colour in $\{\alpha,\beta\}\setminus\{\varphi(v)\}$ appears an odd number of times in $N_{G}(w)$.  In particular, at least one of $\alpha$ or $\beta$ appears an odd number of times in $N_G(w)$.
      \end{itemize}
    \end{itemize}
  \end{itemize}
  Therefore, by choosing $\varphi(v)\not\in X\cup Y$ we obtain an odd colouring of $G$ that satisifies (\ding{96}).  Since $|X\cup Y|\le |X|+|Y|\le 8t+3$, there exists some $\varphi(v)\in \{1,\ldots,8t+4\}\setminus(X\cup Y)$ that completes the colouring of $G$.
\end{proof}

\section{Remarks}

Our proof of \cref{new_result} is inspired by the proof of the result on $d$-degenerate minor-closed families of \citet{cranston.lafferty.ea:note}, which is an inductive proof that involves contracting an edge $wv$ incident to a vertex $v$ of degree at most $d$.  The contraction of this edge (as opposed to the deletion of $v$) is crucial to ensuring that $N_G(v)$ has a colour (namely $\varphi(w)$) that appears an odd number of times.  However, since the class of graphs with product structure is not minor-closed we can not use edge contractions. Instead, we use vertex deletion along with condition (\ding{96}) to achieve a similar effect.

One might hope that \cref{new_result} could be generalized to the setting in which $H$ belongs to some $t$-degenerate minor-closed family of graphs.  However, bounding the size of the set $X$ required to maintain (\ding{96}) when choosing $\varphi(v)$ relies critically on (\ding{74}), which is a property of graphs of treewidth $t$ that is not true for all $t$-degenerate minor-closed graph families.

\paragraph{Subgraphs of $H\boxtimes P\boxtimes K_\ell$}

A number of product structure theorems characterize graphs as subgraphs of $H\boxtimes P\boxtimes K_{\ell}$ where $H$ is a $t$-tree, $P$ is a path, and $K_\ell$ is a complete graph of order $\ell$.  Since $H\boxtimes P\boxtimes K_{\ell}$ is isomorphic to $H\boxtimes K_\ell\boxtimes P$ and $H\boxtimes K_\ell$ is a $(\ell(t+1)-1)$-tree, \cref{new_result} immediately implies that these graphs have odd colourings using $8\ell t+8\ell-4$ colours.

It is possible to improve this slightly by redoing the proof \cref{new_result}.  In this case, the vertex $v:=(x_i,y_j,z_k)$ that is removed is also chose to maximize $(i,j)$, with ties broken arbitrarily. Then one finds that the sizes of the colour sets $X$ and $Y$ that must be avoided when choosing $\varphi(v)$ are bounded by $|X|\le 5\ell t + 3\ell-1$ and $|Y|\le 3\ell t + 2\ell -1$ so that $|X\cup Y|\le 8\ell t + 5\ell -2$. This gives the following variant of \cref{new_result}:

\begin{thm}\label{new_result_kl}
  Let $H$ be a $t$-tree, let $P$ be a path, let $K_\ell$ be a clique on $\ell$ vertices, and let $G$ be a subgraph of $H\boxtimes P\boxtimes K_\ell$.  Then $G$ has an odd colouring using at most $8\ell t + 5\ell -1$ colours.
\end{thm}

\paragraph{Subgraphs of $H\boxtimes I$}

Perhaps a more interesting generalization comes by replacing $P$ by some graph $I$ of maximum-degree $\Delta$.  Again, one can follow the same general strategy used in the proof of \cref{new_result}, with the following changes.
\begin{compactitem}
  \item The vertices $y_1,\ldots,y_h$ are the vertices of $I$ is no particular order.
  \item The set $C_{(x_i,y_i)}$ is defined as
\[
    C_{(x_i,y_i)}:=(\{x_i\}\cup C_{x_i})\boxtimes (\{y_j\}\cup N_I(y_j)) \enspace .
\]
  \item $|X\cup Y|$ is bounded as follows: $|Y|\le|N_G(v)|\le |C_v\setminus\{v\}| \le (t+1)(\Delta+1)-1$.  The set $R$ used to define $X$ is given by $R:=(\{x_i\}\cup C_{x_i})\times N^2_I(y_i)$, where $N^2_I(y_i)$ denotes the set of at most $\Delta^2+1$ vertices in $I$ of distance at most $2$ from $y_i$.  Then $|X|\le|R\setminus\{v\}|\le (t+1)(\Delta^2+1)-1$.  Therefore $|X\cup Y|\le |X|+|Y|\le (\Delta^2+\Delta)(t+1)+2t$.
\end{compactitem}
These changes prove the following variantion of \cref{new_result}:

\begin{thm}\label{new_result_delta}
  Let $H$ be a $t$-tree, let $I$ be a graph of maximum-degree $\Delta$, and let $G$ be a subgraph of $H\boxtimes I$.  The $G$ has an odd colouring using at most $(\Delta^2+\Delta)(t+1)+2t+1$ colours.
\end{thm}

%
% We conclude by noting that the proof of \cref{new_result} can easily be modified to prove the following refinement:
%
% \todo[inline]{We have to rethink this next part now.}
% \begin{thm}\label{new_result2}
%   Let $H$ be a $t$-tree, let $B$ be a $d$-degenerate\footnote{A graph $G$ is $d$-degenerate if every subgraph of $G$ contains a vertex of degree at most $d$.} graph of maximum-degree $\Delta$, and let $G$ be a subgraph of $H\boxtimes B$. Then $G$ has a proper odd colouring $\varphi:V(H)\to\{1,\ldots,2((\Delta+1)t + d)+1\}$.
% \end{thm}
%
% The proof of \cref{new_result2} is almost identical to the proof of \cref{new_result} except for the following modifications:  Let $y_1,\ldots,y_h$ be an ordering of $V(B)$ for which $|N_B(y_j)\cap\{y_1,\ldots,y_{j}\}|\le d$ for each $j\in\{1,\ldots,h\}$ and let $C_{y_j}:=N_B(y_j)\cap\{y_1,\ldots,y_{j}\}$.  Then the definition of $C_{x_i,y_j}$ becomes
% \[
%    C_{x_i,y_j} := V(G)\cap \left( (C_{x_i}\times (\{y_j\}\cup N_B(y_j)))
%        \cup \{(x_i,y): y \in C_{y_j}\}  \right) \enspace .
% \]
% As before the vertex $v:=(x_i,y_j)\in V(G)$ is chosen so that $(i,j)$ is lexicographically maximum.  With this choice of $v$,  $N_G(v)\subseteq C_v$, so $|N_G(v)|\le (\Delta+1)t+d$. Since $|X|\le |C_v|$ and $|Y|\le |N_G(v)|$, this means there are at most $2((\Delta+1)t+d)$ colours that must be avoided when choosing $\varphi(v)$, so $2((\Delta+1)t+d)+1$ colours are sufficient. \Cref{new_result} is a consequence of \cref{new_result2} because a path $P$ is a graph with maximum degree $\Delta=2$ and degeneracy $d=1$.

\bibliographystyle{plainurlnat}
\bibliography{odd}

\end{document}